\documentclass[12pt]{article}

\usepackage{amssymb,amsmath,latexsym}
\usepackage[utf8]{inputenc}
\usepackage{textcomp}
\usepackage[T1]{fontenc}
\usepackage{authblk}
\usepackage{amssymb,amscd,amsmath,amsfonts}
\usepackage{latexsym}
\usepackage{amsthm}
\usepackage{enumerate}
\usepackage{verbatim}
\usepackage{graphicx}
\usepackage[dvips]{epsfig}
\usepackage{subcaption}
\usepackage{amsbsy}
\usepackage{mathtools}
\allowdisplaybreaks

\newtheorem{thm}{Theorem}[section]
\newtheorem{lem}[thm]{Lemma}
\newtheorem{prop}[thm]{Proposition}

\newtheorem{dfn}[thm]{Definition}
\newtheorem{exa}[thm]{Example}
\newtheorem{cor}[thm]{Corollary}
\newtheorem{rk}[thm]{Remark}

\numberwithin{equation}{section}

\newcommand{\C}{\mathbb C}
\newcommand{\R}{\mathbb R}

\newcommand{\N}{\mathbb N}

\newcommand{\e}{\varepsilon}

\newcommand{\bsg}{\boldsymbol{g}}

\newcommand{\tr}{\operatorname{tr}}

\begin{document}

\title{Hyers-Ulam stability of the first order difference equation generated by linear maps}

\author
{Young Woo Nam 
}

\affil[]{\small  
      College of Engineering,
      Konkuk University, 05029 
      Seoul, Korea \\ \texttt{namyoungwoo@konkuk.ac.kr}
      }
\date{}

\maketitle

\begin{abstract}
Hyers-Ulam stability of the difference equation $ z_{n+1} = a_nz_n + b_n $ is investigated. If $ \prod_{j=1}^{n}|a_j| $ has subexponential growth rate, then difference equation generated by linear maps has no Hyers-Ulam stability. Other complementary results are also found where $ \lim_{n \rightarrow \infty}\left( \prod_{j=1}^{n}|a_j| \right)^{\frac{1}{n}} $ is greater or less than one. These results contain Hyers-Ulam stability of the first order linear difference equation with periodic coefficients also. 
\end{abstract}

\section{Introduction}
Ulam suggested in 1940 the existence of an approximate homomorphism between metric groups in \cite{ulam}. Hyers answered this question for the additive functional equation in \cite{hyers} . In decades Hyers-Ulam stability is studied by many researchers and there are fruitful results. More recently, Hyers-Ulam stability has been applied to discrete type equations. Hyers-Ulam stability linear difference equation has been explored by researchers. For example, see \cite{BBP,Po}. Hyers-Ulam stability of $ h- $difference equations are studied in \cite{AO1,AO2,AOR}. Hyers-Ulam stability of a certain type of linear fractional map is searched in \cite{nam3}. Hyers-Ulam stability of functional equation or differential equation is searched in \cite{BCL,HJL}. Eigenvalues determine Hyers-Ulam stability of a certain type of matrices with periodic coefficient in \cite{BOS}. Furthermore, in \cite{BD}, the shadowing properties and Hyers-Ulam stability of linear sequences in Banach spaces are dealt with under the assumption of exponential dichotomy or trichotomy. See \cite{BD} and references therein. \\
\\
In this paper, we discuss Hyers-Ulam stability of the sequence generated by linear maps $ z_{n+1} = a_nz_n + b_n $ 
on complex plane with $ a_n \neq 0 $ for all $ n \in \N $. 
Let $ \{t_n\}_{n\in \N} $ be a sequence of positive real numbers which satisfies that 
$$ \lim_{n \rightarrow \infty} \frac{t_n}{\sum_{j=1}^{n-1}t_j} = 0 . $$
We call that the sequence $ \{t_n\}_{n\in \N} $ has {\em subexponential growth rate}. This sequence is strictly slower than that of any sequence which has exponential expansion or contraction as $ n \rightarrow \infty $. 
We provide a sufficient condition of subexponential growth rate of the linear recurrence, $ z_{n+1} = a_nz_n + b_n $, 
$$
 \lim_{n \rightarrow \infty} \left(\prod_{j=1}^{n}|a_j| \right)^{\frac{1}{n}} = 1 .
$$
If the above limit exists and it is away from one, then the linear recurrence has eventually either at least exponential increasing or decreasing rate. In Section 6, we show that difference equations under this assumption have Hyers-Ulam stability. \\
In this paper we make the contribution to Hyers-Ulam stability of the sequence generated by linear maps. Subexponential growth rate is suggested and a good sufficient condition is provided, which is the limit of the geometric average of $ |a_n| $, the absolute value of coefficients of linear maps. The quantity $ \lim_{n \rightarrow \infty} \left(\prod_{j=1}^{n}|a_j| \right)^{\frac{1}{n}} $ may be a criterion to determine Hyers-Ulam stability of linear recurrence in Section 5 and 6. \\
This paper is organized as follows. In Section 3, a sufficient condition for no Hyers-Ulam stability is suggested with infinite product of the coefficients of linear recurrences. In Section 4, we examine a sufficient condition for Hyers-Ulam stability of linear recurrence with periodic coefficients. In Section 5 and Section 6, we show how the limit of $ \left(\prod_{j=1}^{n}|a_j| \right)^{\frac{1}{n}} $ affects Hyers-Ulam stability of difference equation generated by linear maps.

\section{Preliminaries}
The following Stolz-Ces\`aro theorem \cite{toth} is about the inequality of limit.
\begin{lem}[Stolz-Ces\`aro theorem]
%
Let $ \{x_n\}_{n\in\N} $ and $ \{y_n\}_{n\in\N} $ be real sequences such that $ \{y_n\}_{n\in\N} $ is strictly increasing with $ \lim_{n \rightarrow \infty} y_n = \infty $. Then we have
$$
\liminf_{n \rightarrow \infty} \frac{x_n - x_{n-1}}{y_n-y_{n-1}} \leq \liminf_{n \rightarrow \infty} \frac{x_n}{y_n} \leq \limsup_{n \rightarrow \infty} \frac{x_n}{y_n} \leq \limsup_{n \rightarrow \infty} \frac{x_n - x_{n-1}}{y_n-y_{n-1}} .
$$
In particular, 
$$ \lim_{n \rightarrow \infty} \frac{x_n - x_{n-1}}{y_n-y_{n-1}} = \lim_{n \rightarrow \infty} \frac{x_n}{y_n} $$
provided that the limit of the left side exists.
\end{lem}
%
%
\noindent Let $ y_n = n $ and replace $ x_n - x_{n-1} $ by $ x_n $ in the Stolz-Ces\`aro theorem. Thus we obtain that 
$$
\liminf_{n \rightarrow \infty} x_n \leq \liminf_{n \rightarrow \infty} \frac{x_1 + x_2 + \cdots + x_n}{n} \leq \limsup_{n \rightarrow \infty} \frac{x_1 + x_2 + \cdots + x_n}{n} \leq \limsup_{n \rightarrow \infty} x_n .
$$
\\
\noindent Let $ \{x_n\}_{n\in\N} $ be the sequence of positive real numbers and let $ y_n = \log_2(x_n) $. Thus $ 2^{\frac{y_1 + y_2 + \cdots + y_n}{n}} = \sqrt[n]{x_1\cdots x_n} $. Using the monotonicity of exponentiation and Stolz-Ces\`aro theorem, we have 
$$ \liminf_{n \rightarrow \infty} \frac{x_n}{x_{n-1}} \leq \liminf_{n \rightarrow \infty} \sqrt[n]{x_n} \leq \limsup_{n \rightarrow \infty} \sqrt[n]{x_n} \leq \limsup_{n \rightarrow \infty} \frac{x_n}{x_{n-1}} . $$ 
%
%
\noindent Let $ g(z) $ be the linear map on the complex plane, that is, $ g(z) = az + b $ for $ a,b \in \C $ and $ a \neq 0 $. Then linear maps are classified as follows
\begin{enumerate}
\item $ g $ is parabolic if $ a = \pm 1 $, that is, $ g $ is translation, 
\item $ g $ is elliptic if $ |a| = 1 $ but $ a \neq \pm 1 $, that is, $ g $ is rotation, 
\item $ g $ is loxodromic if $ |a| \neq 1 $,
\item $ g $ is hyperbolic if $ |a| \neq 1 $ and $ a \in \R $, that is, $ g $ is dilation or contraction along straight lines.
\end{enumerate}
For instance, the map $ f(z) = az + b $ is hyperbolic linear map where $ a, b \in \R $ and $ |a| \neq 1 $. This classification is also applied to linear fractional maps in \cite{Beardon}. If $ |a| \neq 1 $, that is, $ g $ is loxodromic, then $ \prod_{j=1}^{n-1}|a| = |a|^{n-1} $, which increases or decreases exponentially by iteration of $ g $. The matrix representation of $ g $ used in Example \ref{eg-no HU stability example1} and this representation $ f(z) = az + b $ of  is as follows
\begin{align*}
\begin{pmatrix}
\sqrt{a} & \frac{b}{\sqrt{a}} \\
0 & \frac{1}{\sqrt{a}} 
\end{pmatrix} \in SL(2,\C)/\pm I 
\end{align*}
where $ SL(2,\C) $ is the special linear group of $ 2 \times 2 $ matrices over $ \C $. Trace of matrix is contained in $ \C \setminus [-2,2] $ if and only if $ g $ is loxodromic, trace is $ \pm 2 $ if and only if $ g $ is parabolic and trace is contained in $ (-2,2) $ if and only if $ g $ is elliptic. 
\\
\noindent We denote the Hyers-Ulam stability by HU stability in the title of section. We use the definition of notions for product and sum as follows for calculations $ \prod_{j=1}^0 a_j = 1, \quad \sum_{j=1}^0 r_j = 0 $
for any complex numbers $ a_j $ and $ r_j $. 


\section{No HU stability of difference equation from linear maps}
%
\medskip 
%
%
%
%
\noindent Let $ g_n(z_n) = a_nz_n + b_n $ be the map on complex plane $ \C $ for complex numbers $ a_n $, $ b_n $ and $ z_n $ for $ n \in \N $ throughout this paper. 
\begin{dfn}
Let $ \{w_n\}_{n \in \N} $ be the sequence which satisfies the inequality
$$ |w_{n+1} - g_n(w_n)| \leq \e $$
for a given $ \e > 0 $ and for all $ n \in \N $. It is called an approximate solution of the difference equation, $ z_{n+1} = g_n(z_n) $. This equation is called Hyers-Ulam stable if there exists a sequence $ \{z_n\}_{n \in \N} $ as follows
\begin{align*}
z_{n+1} = g_n(z_n) \quad \text{and} \quad |w_n - z_n| \leq G(\e)
\end{align*}
for all $ n \in \N $ and $ G(\e) \rightarrow 0 $ as $ \e \rightarrow 0 $. 
\end{dfn}
\medskip
\begin{lem} \label{lem-general term from linear map}
Let the sequence $ \{z_n\}_{n\in \N} $ satisfies that 
$$ z_{n+1} = a_nz_n + b_n $$
for every $ n \in \N $. Then the $ z_n $ is as follows 
\begin{align*} 
z_n = z_1\prod_{j=1}^{n-1}a_j + \sum_{j=0}^{n-2}b_{n-j-1}\prod_{i=1}^{j}a_{n-j+i-1}
\end{align*}
for $ n \geq 2 $ where $ \displaystyle \prod_{i=1}^0a_{n-j+i-1} = 1 $. 
\end{lem}

\begin{proof}
For $ n=2 $, 
$$ z_2 = z_1\prod_{j=1}^{1}a_j + \sum_{j=0}^{0}b_{2-j-1}\prod_{i=1}^{j}a_{2-j+i-1} = a_1z_1 +b_1 . $$
Suppose that $ \displaystyle z_{n-1} = z_1\prod_{j=0}^{n-2}a_j + \sum_{j=1}^{n-3}b_{n-j-2}\prod_{i=1}^{j}a_{n-j+i-2} $. Thus the equation $ z_n = a_{n-1}z_{n-1} + b_{n-1} $ implies that 
\begin{align*}
 z_n &= a_{n-1}z_{n-1} + b_{n-1} \\
&= a_{n-1} \left( z_1\prod_{j=1}^{n-2}a_j + \sum_{j=0}^{n-3}b_{n-j-2}\prod_{i=1}^{j}a_{n-j+i-2} \right) + b_{n-1} \\
&= z_1\prod_{j=1}^{n-1}a_j + \sum_{j=0}^{n-3}b_{n-j-2}\left( a_{n-1}\prod_{i=1}^{j}a_{n-j+i-2} \right) + b_{n-1} \\
&= z_1\prod_{j=0}^{n-1}a_j + \sum_{j=0}^{n-2}b_{n-j-1}\prod_{i=1}^{j}a_{n-j+i-1} .
\end{align*}
Hence, the proof is complete by induction. 
\end{proof}
\medskip 
\noindent Let $ g_n $ be the map which is $ g_n(z) = a_nz+b_n $. Denote $ g_n \circ g_{n-1} \circ \cdots \circ g_1 $ by $ \bsg_n $. Observe that $ \bsg_n = g_n \circ \bsg_{n-1} $. Thus the sequence $ \{z_n\}_{n\in \N} $ in Lemma \ref{lem-general term from linear map} satisfies that $ z_{n+1} = g_n(z_n) $ and $ z_{n+1} = \bsg_n(z_1) $. 
\medskip

\begin{lem} \label{lem-Rn the difference}
Let $ g_n $ be the map which is $ g_n(z) = a_nz+b_n $. Denote the composition of maps, $ g_n \circ g_{n-1} \circ \cdots \circ g_1 $ by $ \bsg_n $ for $ n \in \N $. Then the sequence $ \{w_n\}_{n\in \N} $ is defined by 
$$ w_{n+1} - \bsg_n(w_1) = R_n $$
where $ R_0 = 0 $ for $ n \in \N $ if and only if 
$$ w_{n+1} - g_n(w_n) = R_n - a_nR_{n-1} $$ for $ n \in \N $. 
\end{lem}

\begin{proof}
Assume that $ w_{n+1} - \bsg_n(w_1) = R_n $ for $ n \in \N $. Then the straightforward calculation completes the proof as follows 
\begin{align*}
w_{n+1} - g_n(w_n) &= \bsg_n(w_1) + R_n - (a_nw_n + b_n) \\
&= \bsg_n(w_1) + R_n - \big( a_n(\bsg_{n-1}(w_1) + R_{n-1}) +b_n \big) \\
&= \bsg_n(w_1) - \big( a_n \bsg_{n-1}(w_1) +b_n \big) + R_n - a_nR_{n-1} \\
&= R_n - a_nR_{n-1} .
\end{align*}
Conversely, for $ n=1 $ we have $ w_2 - g_1(w_1) = R_1 - a_1R_0 = R_1 $. Assume for induction that $ w_n - \bsg_{n-1}(w_1) = R_{n-1} $. Recall that $ \bsg_n = g_n \circ \bsg_{n-1} $ and $ g_n(z) = a_nz + b_n $. Then we have that
\begin{align*}
w_{n+1} - \bsg_n(w_1) &= g_n(w_n) + R_n - a_nR_{n-1} - g_n(\bsg_{n-1}(w_1)) \\
&= g_n(w_n) + R_n - a_nR_{n-1} - g_n(w_n -R_{n-1}) \\
&= a_nw_n + b_n + R_n - a_nR_{n-1} - \big( a_n(w_n-R_{n-1}) + b_n \big) \\
&= R_n .
\end{align*}
The proof is complete by induction. 
\end{proof}
\medskip
%
%
\begin{lem} \label{lem-difference between w and z}
Let $ g_n(z) =a_nz + b_n $ be the map on complex plane for $ n \in \N $. Let $ \{z_n\}_{n\in \N} $ and $ \{w_n\}_{n\in \N} $ be the sequences defined as follows 
$$ z_{n+1} = g_n(z_n) 
, \qquad w_{n+1} = g_n(w_n) + r_n 
$$
respectively for $ n \in \N $. Then the following equation holds
\begin{align*}
w_{n+1} -z_{n+1} 
&= \sum_{j=0}^{n-1} r_{n-j} \prod_{i=1}^j a_{n-j+i} + (w_1 -z_1) \prod_{j=1}^n a_j .
\end{align*}
\end{lem}

\begin{proof}
Lemma \ref{lem-Rn the difference} implies that $ w_{n+1} - \bsg_n(w_1)= R_n $ where $ r_n = R_n - a_nR_{n-1} $ and $ R_0 = 0 $ for $ n \in \N $. The equation $ z_{n+1} = g_n(z_n) = \bsg_n(z_1) $ holds for $ n \in \N $. Thus Lemma \ref{lem-general term from linear map} and Lemma \ref{lem-Rn the difference} yield the following equation
\begin{align*}
w_{n+1} -z_{n+1} &= w_{n+1} - \bsg_n(z_1) \\
&= w_{n+1} - \bsg_n(w_1) + \bsg_n(w_1) - \bsg_n(z_1) \\
&= R_n + (w_1 - z_1) \prod_{j=1}^n a_j 
\end{align*}
for $ n \in \N $. The equation $ r_n = R_n - a_nR_{n-1} $ is similar to the recurrence $ z_{n+1} = a_nz_n + b_n $. Observe that $ R_1 = r_1 $. Then we obtain the following equation by induction 
\begin{align} \label{eq-the equation of Rn}
R_n &= r_n + a_nr_{n-1} + a_na_{n-1}r_{n-2} + \cdots + a_na_{n-1} \ldots a_2r_1 \nonumber \\
&= \sum_{j=0}^{n-1} r_{n-j} \prod_{i=1}^j a_{n-j+i}
\end{align}
where $ \displaystyle \prod_{i=1}^0 a_{n-j+i} = 1 $. The proof is complete.
\end{proof}
%
 
\medskip

\begin{thm} \label{thm-non stability of linear maps}
Let $ \{z_n\}_{n\in \N} $ be the sequence defined as follows
\begin{align} \label{eq-equation of z}
z_{n+1} = a_nz_n + b_n
\end{align}
for $ n \in \N $. Suppose that the following conditions  
$$ 
 \sup_{n \in \N} \left\{\!\,\prod_{j=1}^{n-1} |a_j| \!\,\right\} < \infty \quad \text{and} \quad  \sup_{n \in \N} \left\{\!\,n  \prod_{j=1}^{n-1} |a_j| \!\,\right\} = \infty. 
$$ 
Then the equation \eqref{eq-equation of z} has no Hyers-Ulam stability. 
\end{thm}

\begin{proof}
Let $ \{w_n\}_{n\in \N} $ be the sequence defined by the equation
$$ w_{n+1} - g_n(w_n) = r_n $$
where $ g_n(z) = a_nz +b_n $. For any given $ \e>0 $, choose $ r_n $ which satisfies the following equation 
\begin{align} \label{eq-rn defined for theorem in sec3}
r_n = \frac{\e}{M} \prod_{j=1}^{n} a_j \ \ \text{where} \ \ M = \sup_{n \in \N} \left\{\!\,\prod_{j=1}^{n-1} |a_j| \!\,\right\} .
\end{align}
%
Thus $ |r_n| \leq \e $ for all $ n \in \N $. Define $ R_n $ using the recursive formula $ R_n = a_nR_{n-1} + r_n $ and $ R_0 = 0 $ for $ n \in \N $. Recall that $ \bsg_n = g_n \circ g_{n-1} \circ \cdots \circ g_1 $ for $ n \in \N $. Thus Lemma \ref{lem-Rn the difference} implies that 
$$ 
w_{n+1} - \bsg_n(w_1) = R_n .
$$
Moreover, the equation \eqref{eq-the equation of Rn} and defined $ r_n $ in \eqref{eq-rn defined for theorem in sec3} imply that
\begin{align*}
R_n 
&= \sum_{j=0}^{n-1} r_{n-j} \prod_{i=1}^j a_{n-j+i}
= \sum_{j=0}^{n-1} \frac{\e}{M} \prod_{j=1}^{n-j} a_j \prod_{i=1}^j a_{n-j+i} 
= \frac{\e}{M} \sum_{k=0}^{n-1} \prod_{j=1}^{n} a_j 
= \frac{\e n}{M} \prod_{j=1}^{n} a_j .
\end{align*}
Lemma \ref{lem-difference between w and z} implies that 
\begin{align*}
w_{n+1} - z_{n+1} &= R_n + (w_1 - z_1) \prod_{j=1}^n a_j = \frac{\e}{M}n \prod_{j=1}^{n} a_j + (w_1 - z_1) \prod_{j=1}^n a_j .
\end{align*}
Then $ | w_{n+1} - z_{n+1} | $ is unbounded for $ n \in \N $. Hence,  the equation \eqref{eq-equation of z} has no Hyers-Ulam stability. 
\end{proof}

\medskip

\begin{cor} \label{cor-non stability of linear maps}
Let $ \{z_n\}_{n\in \N} $ be the sequence defined as follows
$$ z_{n+1} = a_nz_n + b_n $$
for $ n \in \N $. Suppose that the inequality holds as follows 
$$ 0 < \inf_{n \in \N} \prod_{j=1}^{n-1} |a_j| \leq \sup_{n \in \N} \prod_{j=1}^{n-1} |a_j| < \infty . $$ 
Then the above recurrence $ z_{n+1} = a_nz_n + b_n $ has no Hyers-Ulam stability. 
\end{cor}

\begin{proof}
The inequality $ 0 < \inf_{n \in \N} \prod_{j=1}^{n-1} |a_j| $ implies that 
$$ \lim_{n \rightarrow \infty} n \prod_{j=1}^{n-1} |a_j| = \infty . $$
\end{proof}


\begin{exa} \label{eg-no HU stability example1}
Let $ g_n $ be the map as follows 
\begin{align} \label{eq-map for example1 in sec3}
g_n(z) = \left(1+\frac{1}{n^2} \right)^2 e^{2\pi \alpha i}(z-2)
\end{align}
for the real number $ \alpha \in [0,1) $. Denote $ \left(1+\frac{1}{n^2} \right)^2e^{2\pi \alpha i} $ by $ a_n $. Thus we have
$$
\sum_{n=1}^{\infty} \big| |a_n|-1 \big| \leq \sum_{n=1}^{\infty} \frac{3}{n^2} < \infty .
$$
The constant $ \prod_{n=1}^{\infty}|a_n| $ is also convergent and has non zero limit. %
Let $ \{z_n\}_{n\in \N} $ be the sequence satisfying that 
$$
z_{n+1} = g_n(z_n)
$$
for $ n \in \N $. Then Theorem \ref{thm-non stability of linear maps} implies that the difference equation defined by the equation \eqref{eq-map for example1 in sec3} has no Hyers-Ulam stability. 
The matrix representation of $ g_n $ is as follows 
\begin{align*}
\begin{pmatrix}
\frac{n^2+1}{n^2}e^{\pi \alpha i} & \frac{-2(n^2+1)}{n^2}e^{\pi \alpha i} \\[0.5em]
0 & \frac{n^2}{n^2+1}  e^{-\pi \alpha i} 
\end{pmatrix} .
\end{align*}
Thus the trace of the above matrix is that $ \tr(g_n) = \frac{n^2+1}{n^2}e^{\pi \alpha i} + \frac{n^2}{n^2+1}  e^{-\pi \alpha i} $. The limit of $ \tr(g_n) $ is as follows 
$$ \lim_{n \rightarrow \infty} \tr(g_n) = e^{\pi \alpha i} + e^{-\pi\alpha i} = 2\cos(\alpha\pi) . $$
Then $ \lim_{n \rightarrow \infty} \tr(g_n) \in [-2,2] $. 
The maps have two different cases.
\begin{enumerate}
\item If $ \alpha = 0 $, then $ g_n $ converges to the parabolic linear map $ z \mapsto z-2 $ as $ n \rightarrow \infty $. 
\item If $ \alpha \neq 0 $, then $ g_n $ converges to the elliptic linear map $ z \mapsto e^{\pi \alpha i}(z-2) $ as $ n \rightarrow \infty $. 
\end{enumerate}
\end{exa}

\noindent In the above example, every $ g_n $ is the loxodromic linear map. However, if $ g_n $ is arbitrary close to parabolic or elliptic linear map for all $ n \geq N $, then the difference equation $ z_{n+1} = g_n(z_n) $ from the map \eqref{eq-map for example1 in sec3} does not have Hyers-Ulam stability. 
\medskip

\begin{exa} \label{exa-loxo generated but no stable}
Let $ g_n $ be the map as follows 
\begin{align} \label{eq-map for example2 in sec3}
g_n(z) = a_nz + 5,\quad \text{where} \quad a_n =
\begin{cases}
2, \quad n = 2k - 1 \\
\frac{1}{2}, \quad n = 2k
\end{cases}
\end{align}
for $ k \in \N $. Thus we obtain that 
$$ \prod_{j=1}^{2k} a_j = 1 \quad \text{and} \quad \prod_{j=1}^{2k-1} a_j = 2 $$
for every $ k \in \N $. Then the following inequality holds
$$ 
0 < \inf_{n \in \N} \prod_{j=1}^{n-1} |a_j| < \sup_{n \in \N} \prod_{j=1}^{n-1} |a_j| < \infty . 
$$
Let $ \{z_n\}_{n\in \N} $ be the sequence satisfying that 
$$
z_{n+1} = g_n(z_n)
$$
for $ n \in \N $. For any positive $ \e $, choose the sequence $ \{w_n\}_{n \in \N} $
which satisfies the difference equation $ w_{n+1} - g_n(w_n) = \frac{\e}{2}\prod_{j=1}^{n-1}|a_j| $. Then Corollary \ref{cor-non stability of linear maps} implies that the difference equation $ z_{n+1} = g_n(z_n) $ from the map \eqref{eq-map for example2 in sec3} has no Hyers-Ulam stability. 
\end{exa}
\noindent In Example \ref{exa-loxo generated but no stable}, the sequence $ \{z_n\}_{n \in \N} $ with the initial point $ z_0 \in \C $ has the subsequence $ \{z_{3n}\}_{n \in \N} $ which satisfies that 
$$ z_{2k+2} = 2z_{2k-1} + 20, \quad z_{2k+3} = \frac{1}{2}z_{2k} + \frac{25}{2} $$
for $ k \in \N $. Thus each of above subsequences is generated by single hyperbolic linear map and it has Hyers-Ulam stability because of Theorem \ref{thm-stability of exponential shrinking} and Theorem \ref{thm-stability of exponential expanding}. Similarly, difference equations from the subsequences $ \{z_{5n}\}_{n \in \N} $ and $ \{z_{7n}\}_{n \in \N} $ also have Hyers-Ulam stability.

\section{HU stability of difference equation from linear maps}
%

\begin{prop} \label{prop-stability of linear maps}
Let $ \{z_n\}_{n\in \N} $ be the sequence defined as follows
\begin{align} \label{eq-linear map in sec4}
z_{n+1} = a_nz_n + b_n
\end{align}
for $ n \in \N $. Suppose that the sum 
$$ 
1 + |a_n| + |a_na_{n-1}| + \cdots + |a_na_{n-1}\ldots a_2|  
$$
is bounded for all $ n \in \N $. Then the difference equation \eqref{eq-linear map in sec4} has Hyers-Ulam stability.
\end{prop}

\begin{proof}
Let $ g_n(z) = a_nz + b_n $ for every $ n \in \N $. Thus $ z_{n+1} = g_n(z_n) $ for every $ n \in \N $. Let $ \{w_n\}_{n \in \N} $ be the sequence which satisfies that the recursive relation 
$$ w_{n+1} - g_n(w_n) = r_n $$
for every $ n \in \N $. Lemma \ref{lem-difference between w and z} implies that 
\begin{align} \label{eq-difference between w and z again}
w_{n+1} -z_{n+1} &= \sum_{j=0}^{n-1} r_{n-j} \prod_{i=1}^j a_{n-j+i} + (w_1 -z_1) \prod_{j=1}^n a_j 
\end{align}
for $ n \in \N $. Thus the equation \eqref{eq-the equation of Rn} yields that
\begin{align} 
\sum_{j=0}^{n-1} r_{n-j} \prod_{i=1}^j a_{n-j+i} = r_n + a_nr_{n-1} + a_na_{n-1}r_{n-2} + \cdots + a_na_{n-1} \ldots a_2r_1 . \nonumber 
\end{align}
%
For any given $ \e > 0 $, assume that $ |r_n| \leq \e $ for all $ n \in \N $. The equation \eqref{eq-difference between w and z again} implies that
\begin{align*} 
|w_{n+1} -z_{n+1}| &\leq \big(1 +| a_n| + |a_na_{n-1}| + \cdots + |a_na_{n-1}\ldots a_2| \big)\;\!\e \nonumber \\
&\qquad + |w_1 -z_1| \prod_{j=1}^n |a_j| .
\end{align*} 
Then 
we may choose $ z_1 = w_1 $ and $ \{z_n\}_{n\in \N} $ as a solution for \eqref{eq-linear map in sec4}. Hence, the equation \eqref{eq-linear map in sec4} has Hyers-Ulam stability. 
\end{proof}
\medskip
\begin{rk}
In the proof of Proposition \ref{prop-stability of linear maps}, the positive number, $ \prod_{j=1}^n |a_j| $ is bounded for all $ n \in \N $ by assumption. Then we may choose $ z_1 $ close enough to $ w_1 $ depending on $ \e $. Thus the sequence $ \{z_n\}_{n\in \N} $ is not unique for the given approximate solution of \eqref{eq-linear map in sec4}. 
\end{rk}
\noindent If all numbers $ a_n $ are the same as $ a $ and $ 0 < |a| < 1 $, then the sum
$$ 1 + a_n + a_na_{n-1} + \cdots + a_na_{n-1}\ldots a_2 $$
is a geometric series with the common ratio $ a $. Proposition \ref{prop-stability of linear maps} implies that the difference equation $ z_{n+1} = az_n + b_n $ with $ |a| < 1 $ has Hyers-Ulam stability. 
\medskip
\begin{thm} \label{thm-HU stability with periodic coefficients 1}
Let $ \{z_n\}_{n\in \N} $ be the sequence defined as follows
\begin{align} \label{eq-equation of z in sec4 thm}
z_{n+1} = a_nz_n + b_n
\end{align}
for $ n \in \N $. Let $ a_n $ be the periodic coefficient, that is, $ a_{n+p} = a_n $ for some $ p \in \N $ and for every $ n \in \N $. Assume that $ \prod_{k=1}^p \left| a_k \right| < 1 $. Then the difference equation \eqref{eq-equation of z in sec4 thm} has Hyers-Ulam stability.
\end{thm}

\begin{proof}
It suffice to show that $ | 1 + a_n + a_na_{n-1} + \cdots + a_na_{n-1}\ldots a_2 | $ is bounded for all $ n \in \N $ by Proposition \ref{prop-stability of linear maps}. Let $ n = mp + l $ where $ m \in \N_0 $ and $ l = 0,1,2, \ldots , p-1 $. Use the following notions for later calculations
\begin{align*}
q &= \prod_{k=1}^p a_k, \\ 
S &= 1 + a_p + a_pa_{p-1} + \cdots + a_pa_{p-1} \cdots a_2, \\
S_m &= 1 + a_{mp} + a_{mp}a_{mp-1} + \cdots + a_{mp}a_{mp-1} \cdots a_2
\end{align*}
for $ m \in \N $. Observe that $ |q| < 1 $ and $ S_1 = S $. 
Thus periodicity of $ a_n $ implies that 
\begin{align*}
S_m &= 1+ a_{mp} + a_{mp}a_{mp-1} + \cdots + a_{mp}a_{mp-1} \cdots a_2 \\
& = 1+ a_{mp} + a_{mp}a_{mp-1} + \cdots + a_{mp}a_{mp-1} \cdots a_{(m-1)p+2} \\
& \qquad \ \ + a_{mp}a_{mp-1} \cdots a_{(m-1)p+1} \big(1+ a_{(m-1)p} + a_{(m-1)p}a_{(m-1)p-1} + \\
& \qquad \qquad \cdots + a_{(m-1)p}a_{(m-1)p-1} \cdots a_{2}  \big) \\
& = S + q\big(1+ a_{(m-1)p} + a_{(m-1)p}a_{(m-1)p-1} + 
\cdots + a_{(m-1)p}a_{(m-1)p-1} \cdots a_{2}  \big) \\
&= S + qS_{m-1} .
\end{align*}
Then the recursive relation holds $ 
S_m - \frac{S}{1-q} = q \left( S_{m-1} - \frac{S}{1-q} \right) $ and then
\begin{align*} 
S_m = q^{m-1} \left(S - \frac{S}{1-q} \right) + \frac{S}{1-q} = \frac{1-q^m}{1-q}S
\end{align*}
for $ m \in \N $. The sum $ |S_m| $ is bounded above by $  \frac{2}{|1-q|} |S| $ for all $ m \in \N $. Moreover, the periodicity of $ a_n $ also implies that
\begin{align} \label{eq-sum with Sm}
& 1 + a_n + a_na_{n-1} + \cdots + a_na_{n-1}\ldots a_2 \nonumber \\
& \quad = 1 + a_n + a_na_{n-1} + \cdots + a_na_{n-1}\ldots a_{n-l+2} \nonumber \\
& \qquad \ + a_na_{n-1}\ldots a_{n-l+2}a_{n-l+1} \big(1+ a_{mp} + a_{mp}a_{mp-1} + \cdots + a_{mp}a_{mp-1} \cdots a_2 \big) \nonumber \\
& \quad = 1 + a_l + a_la_{l-1} + \cdots + a_la_{l-1} \cdots a_2 \nonumber \\
& \qquad \ + a_la_{l-1} \cdots a_1\big(1+ a_{mp} + a_{mp}a_{mp-1} + \cdots + a_{mp}a_{mp-1} \cdots a_2 \big) \nonumber \\
& \quad = 1 + a_l + a_la_{l-1} + \cdots + a_la_{l-1} \cdots a_3 + a_la_{l-1} \cdots a_2 S_m .
\end{align}
%
Denote the sum of absolute values of each term of $ S_m $ by $ T_m $, that is,
$$ T_m = 1 + |a_{mp}| + |a_{mp}a_{mp-1}| + \cdots + |a_{mp}a_{mp-1} \cdots a_2| $$
for $ m \in \N $ and denote $ T_1 $ by $ T $. Thus by similar calculation for $ S_m $, we obtain the following inequality
\begin{align*} 
T_m = |q|^{m-1} \left(T - \frac{T}{1-|q|} \right) + \frac{T}{1-|q|} = \frac{1-|q|^m}{1-|q|}T \leq \frac{2}{1-|q|}T
\end{align*}
for all $ m \in \N $. Moreover, the equation which is similar to \eqref{eq-sum with Sm} is satisfied as follows
\begin{align*} 
& 1 + |a_n| + |a_na_{n-1}| + \cdots + |a_na_{n-1}\ldots a_2| \nonumber \\
& \quad = 1 + |a_l| + |a_la_{l-1}| + \cdots + |a_la_{l-1} \cdots a_3| + |a_la_{l-1} \cdots a_2| \;\!T_m .
\end{align*}
Since $ 1 \leq l \leq p-1 $ and $ T_m $ is bounded above for all $ m \in \N $, the number $ 1 + |a_n| + |a_na_{n-1}| + \cdots + |a_na_{n-1}\ldots a_2| $ is also bounded. Hence, Proposition \ref{prop-stability of linear maps} implies that the difference equation \eqref{eq-equation of z in sec4 thm} has Hyers-Ulam stability. 
\end{proof}
\medskip
\begin{rk}
The assumption $ \prod_{k=1}^p \left| a_k \right| < 1 $ of Theorem \ref{thm-HU stability with periodic coefficients 1} cannot be improved. See Example \ref{exa-loxo generated but no stable}.
\end{rk}
\begin{exa}
Let $ g_n $ be the map as follows 
$$
g_n(z) = a_nz + 5,\quad \text{where} \quad a_n =
\begin{cases}
2, \quad n = 3k - 2 \\
i, \quad\, n = 3k-1 \\
\frac{1}{3}, \quad n = 3k
\end{cases}
$$
for $ k \in \N $. 
The sequence of coefficients $ \{a_n\}_{n\in \N} $ is periodic with period 3. The coefficients are $ a_1 = 2 $, $ a_2 = i $ and $ a_3 = \frac{1}{3} $. The product of these numbers is $ |a_1a_2a_3| = \frac{2}{3} < 1 $. Let $ \{z_n\}_{n\in \N} $ be the sequence defined by $ z_{n+1} = g_n(z_n) $ for $ n \in \N $. Theorem \ref{thm-HU stability with periodic coefficients 1} implies that the recurrence $ z_{n+1} = g_n(z_n) $ has Hyers-Ulam stability. Let $ \{w_n\}_{n\in \N} $ be the sequence as follows
$$ w_{n+1} = g_n(w_n) + r_n $$
where $ |r_n| \leq \e $ for all $ n \in \N $. Choose $ z_1 = w_1 $. Thus Proposition \ref{prop-stability of linear maps} and Theorem \ref{thm-HU stability with periodic coefficients 1} imply that 
\begin{align*}
|w_{n+1} - z_{n+1}| &\leq \big(1 +| a_n| + |a_na_{n-1}| + \cdots + |a_na_{n-1}\ldots a_2| \big)\;\!\e \\
 &\leq \big( 1+ |a_2|\;\!T_m \big)\;\!\e \\
 &\leq \big( 1+ 6 \big)\;\!\e \\
 &= 7\;\!\e
\end{align*}
for all $ n \in \N $. Observe that $ T = 1 + |a_3| + |a_3a_2| = 1 + \frac{1}{3} + 1 = \frac{7}{3} $ and $ T_m \leq  \frac{2T}{1-\frac{2}{3}} = 6 $ for all $ m \in \N $. 
\end{exa}

\section{Subexponential growth rate and HU stability} \label{sec-subexponential}
Let $ \{t_n\}_{n \in \N} $ be a sequence of positive real numbers. We suggest that if $ \{t_n\}_{n \in \N} $ satisfies the following limit
$ \lim_{n \rightarrow \infty} \frac{t_n}{\sum_{j=1}^{n-1}t_j} = 0  $, then we call that the sequence $ \{t_n\}_{n \in \N} $ has {\em subexponential growth rate}. 

\medskip 

\begin{lem}[Key lemma 1] \label{lem-key lemma1}
Let $ \{t_n\}_{n \in \N} $ be a sequence of positive real numbers. Suppose that 
$$ \lim_{n \rightarrow \infty} t_n^{\frac{1}{n}} = 1 . $$
Then $ \{t_n\}_{n \in \N} $ has sub-exponential growth rate. Equivalently, the following limit holds
$$ \lim_{n \rightarrow \infty} \frac{t_n}{\sum\limits_{j=1}^{n-1}t_j} = 0 . $$

\end{lem}
\begin{proof}
The assumption is equivalent with
$$ \lim_{n \rightarrow \infty} \frac{\log t_n}{n} = 0 . $$
We claim that for all $ \delta > 0 $ there exists $ n = n(\delta) $ which satisfies that
\begin{align} \label{eq-inequality of delta}
\sum_{j=1}^n t_j < e^{\delta}\sum_{j=1}^{n-1} t_j .
\end{align}
If not, there exists $ \delta > 0 $ which satisfies that
\begin{align} \label{eq-inequality of delta reverse}
\sum_{j=1}^n t_j \geq e^{\delta}\sum_{j=1}^{n-1} t_j .
\end{align}
for all $ n \in \N $, that is, $ t_n \geq (e^{\delta} -1)\sum_{j=1}^{n-1} t_j $. We claim that 
\begin{align} \label{eq-inequality by induction 1}
t_n \geq (e^{\delta} -1)e^{(n-2)\;\!\delta}t_1
\end{align}
for $ n = 2,3,4, \ldots $ by induction. For $ n=2 $, $ t_2 \geq (e^{\delta}-1)t_1 $ by the inequality \eqref{eq-inequality of delta reverse}. Assume that $ t_k \geq (e^{\delta} -1)e^{(k-2)\;\!\delta}t_1 $ for $ k = 2,3,4, \ldots, n $. Thus we obtain the inequality for $ t_{n+1} $ as follows
\begin{align*}
t_{n+1} &\geq e^{\delta}\sum_{j=1}^{n} t_j \\
&= (e^{\delta}-1)(t_1 + t_2 + \cdots + t_n) \\
&\geq (e^{\delta}-1)\big(t_1+ (e^{\delta}-1)t_1 + (e^{\delta}-1)e^{\delta}t_1 + \cdots + (e^{\delta}-1)e^{(n-2)\delta}t_1 \big) \\
&= (e^{\delta}-1)\big(1+ (e^{\delta}-1) + (e^{\delta}-1)e^{\delta} + \cdots + (e^{\delta}-1)e^{(n-2)\delta} \big)t_1 \\
&= (e^{\delta}-1)\big(1+ (-1+e^{\delta}) + (-e^\delta+e^{2\delta}) + \cdots + (-e^{(n-2)\delta} + e^{(n-1)\delta}) \big)t_1 \\
&= (e^{\delta}-1)e^{(n-1)\delta}t_1 .
\end{align*}
Thus $ t_{n+1} \geq (e^{\delta}-1)e^{(n-1)\delta}t_1 $. The claim is proved. Then the equation \eqref{eq-inequality by induction 1} implies the following inequality 
$$ \lim_{n \rightarrow \infty} \frac{\log t_n}{n} \geq {\delta} > 0 , $$
which is the contradiction. The equation \eqref{eq-inequality of delta} implies that there exists $ n $ such that 
$$ 1 < \frac{\sum_{j=1}^n t_j}{\sum_{j=1}^{n-1} t_j} < e^{\delta} $$
for all $ \delta > 0 $. Choose $ \delta = \delta(n) \searrow 0 $ as $ n \rightarrow \infty $. Passing the limit we obtain that $ \lim_{n \rightarrow \infty} \frac{\sum_{j=1}^n t_j}{\sum_{j=1}^{n-1} t_j} = 1 $. Hence, the following limit holds
$$ \lim_{n \rightarrow \infty} \frac{\sum_{j=1}^n t_j}{\sum_{j=1}^{n-1} t_j} - 1 = \lim_{n \rightarrow \infty} \frac{\sum_{j=1}^n t_j  - {\sum_{j=1}^{n-1}t_j}}{\sum_{j=1}^{n-1} t_j} = \lim_{n \rightarrow \infty} \frac{t_n}{\sum_{j=1}^{n-1}t_j} = 0 . $$
\end{proof}

\medskip

\begin{lem} \label{lem-key lemma1 supplementary}
Let $ \{t_n\}_{n \in \N} $ be a sequence of positive real numbers. Suppose that 
$$ \lim_{n \rightarrow \infty} t_n^{\frac{1}{n}} = 1 $$
and suppose also that $ \displaystyle \lim_{n \rightarrow \infty} \sum_{j=1}^n t_j $ exists. 
Then the following limit is satisfied
$$ \lim_{n \rightarrow \infty} \frac{t_n}{\sum\limits_{j=n+1}^{\infty}t_j} = 0 . $$
\end{lem}

\begin{proof}
The equation
$$ \frac{\sum_{j=n}^{\infty}t_j}{\sum_{j=n+1}^{\infty}t_j} = 1 + \frac{t_n}{\sum_{j=n+1}^{\infty}t_j} $$
implies that it suffice to show that 
$$ \frac{\sum_{j=n}^{\infty}t_j}{\sum_{j=n+1}^{\infty}t_j} = 1 . $$
We claim that $ \displaystyle \sum_{j=n}^{\infty}t_j < e^{\delta} \sum_{j=n+1}^{\infty}t_j $ for all $ \delta > 0 $. If not, there exists $ \delta > 0 $ which satisfies that 
\begin{align} \label{eq-inequality reverse}
\sum_{j=n}^{\infty}t_j \geq e^{\delta} \sum_{j=n+1}^{\infty}t_j
\end{align}
for all big enough $ n \in \N $. Thus $ t_n \geq (e^{\delta} -1)\sum_{j=n+1}^{\infty}t_j $. In other words, we obtain that $ t_{m-1} \geq (e^{\delta} -1)\sum_{j=m}^{\infty}t_j $, that is, $ \frac{t_{m-1}}{e^{\delta} -1} \geq \sum_{j=m}^{\infty}t_j $. Then the inequality is satisfied as follows
\begin{align*}
\frac{\log t_m}{m} \leq \frac{\log \left( \sum_{j=m}^{\infty}t_j \right)}{m} \leq \frac{\log (t_{m-1}) - \log (e^{\delta} -1)}{m} .
\end{align*}
Passing the limit by squeezing, we obtain the limit as follows
\begin{align} \label{eq-limit of log seq}
\lim_{m \rightarrow \infty} \frac{\log \left( \sum_{j=m}^{\infty}t_j \right)}{m} = 0 .
\end{align}
We claim that
\begin{align} \label{eq-inequality of t-n reverse order}
t_{n} \geq (e^{\delta} -1)e^{k\delta} \sum_{j=n+k+1}^{\infty}t_j
\end{align} 
for $ k = 0,1,2, \ldots $. For $ k = 0 $, the inequality is the same as \eqref{eq-inequality reverse}, that is,
\begin{align*}
t_n \geq (e^{\delta} -1)\sum_{j=n+1}^{\infty}t_j .
\end{align*}
Denote $ \sum_{j=i}^{\infty}t_j $ by $ A_i $ for $ i \in \N $. Observe that $ A_i - A_{i+1} = t_i $. Thus  the equation \eqref{eq-inequality reverse} yields the inequality $ A_{n+m} \geq e^{\delta}A_{n+m+1} $ for $ m = 0,1,2,\ldots $. Then the inequality is satisfied as follows by comparison of each terms of left and right hand sides 
\begin{align*}
 & t_n + (e^{\delta}-1)A_{n+1} + (e^{\delta}-1)e^{\delta}A_{n+2} + \cdots + (e^{\delta}-1)e^{(k-1)\delta}A_{n+k} \\
 & \qquad \geq (e^{\delta}-1)A_{n+1} + (e^{\delta}-1)e^{\delta}A_{n+2} + \cdots + (e^{\delta}-1)e^{(k-1)\delta}A_{n+k} \\
 & \qquad \qquad + (e^{\delta}-1)e^{k\delta}A_{n+k+1}
\end{align*}
Eliminate same terms in both sides of the above inequality. The proof of \eqref{eq-inequality of t-n reverse order} is done. In the inequality \eqref{eq-inequality of t-n reverse order}, if we choose $ k = n $, then 
$$ t_{n} \geq (e^{\delta} -1)e^{n\delta} \sum_{j=2n+1}^{\infty}t_j $$
This inequality yields that
\begin{align*}
\frac{\log t_n}{n} &\geq \frac{\log (e^{\delta} -1) + n\delta + \log \left(\sum_{j=2n+1}^{\infty}t_j \right) }{n} \\
&= \frac{\log (e^{\delta} -1)}{n} + \delta + \frac{\log \left(\sum_{j=2n+1}^{\infty}t_j \right)}{2n + 1}\cdot \frac{2n+1}{n}
\end{align*}
for all big enough $ n \in \N $. Passing the limit using the inequality \eqref{eq-limit of log seq} as $ n \rightarrow \infty $, we obtain that $ \lim_{n \rightarrow \infty} \frac{\log t_n}{n} \geq \delta > 0 $, which is the contradiction. The claim is proved. Hence, 
$$ \sum_{j=n}^{\infty}t_j < e^{\delta} \sum_{j=n+1}^{\infty}t_j, \ \ \text{that is,} \ \ 1 < \frac{\sum_{j=n}^{\infty}t_j}{\sum_{j=n+1}^{\infty}t_j} < e^{\delta}  $$
for all $ \delta > 0 $. Choose $ \delta = \delta(n) \searrow 0 $ as $ n \rightarrow \infty $. Passing the limit we obtain the limit $ \lim_{n \rightarrow \infty} \frac{\sum_{j=n}^{\infty}t_j}{\sum_{j=n+1}^{\infty}t_j} =1 $. Hence, the following limit is satisfied
\begin{align*}
\lim_{n \rightarrow \infty}\frac{\sum_{j=n}^{\infty}t_j}{\sum_{j=n+1}^{\infty}t_j} - 1 = \lim_{n \rightarrow \infty}\frac{t_n}{\sum_{j=n+1}^{\infty}t_j} = 0 .
\end{align*}
\end{proof}


\begin{thm} \label{thm-no HU stability subexp}
Let $ \{z_n\}_{n\in \N} $ be the sequence defined as follows
\begin{align} \label{eq-equation of z in sec5 thm}
z_{n+1} = a_nz_n + b_n
\end{align}
for $ n \in \N $. Assume that
$$ \lim_{n \rightarrow \infty} \left(\prod_{j=1}^{n} |a_j| \right)^{\frac{1}{n}} = 1 . $$
Then the difference equation \eqref{eq-equation of z in sec5 thm} has no Hyers-Ulam stability. 
\end{thm}

\begin{proof}
For a given $ \e > 0 $, let $ \{w_n\}_{n\in \N} $ be the sequence defined as follows
$$ w_{n+1} = a_nw_n + b_n + r_n $$
where $ r_n = \left(\prod_{i=1}^n \frac{a_i}{|a_i|}\right) \e $ for $ n \in \N $. Lemma \ref{lem-difference between w and z} implies that
\begin{align*}
w_{n+1} -z_{n+1} &= \sum_{j=0}^{n-1} r_{n-j} \prod_{i=1}^j a_{n-j+i} + (w_1 -z_1) \prod_{j=1}^n a_j \\
&= \sum_{j=1}^{n} r_{j} \prod_{i=j+1}^n a_{i} + (w_1 -z_1) \prod_{i=1}^n a_i \\
&= \left( \prod_{i=1}^n a_i \right) \left[ \sum_{j=1}^{n} r_{j} \prod_{i=1}^n a_i^{-1} \prod_{i=j+1}^n a_{i} + (w_1 -z_1) \right] \\
&= \left( \prod_{i=1}^n a_i \right) \left[ \sum_{j=1}^{n} r_{j} \prod_{i=1}^j a_i^{-1} + (w_1 -z_1) \right] \\
&= \left( \prod_{i=1}^n a_i \right) \left[ \e\sum_{j=1}^{n} \prod_{i=1}^j |a_i|^{-1} + (w_1 -z_1) \right]  .
\end{align*}
Thus we obtain that 
\begin{align} \label{eq-inquality of norm wn-zn}
|w_{n+1} -z_{n+1}|
& = \left( \prod_{i=1}^n |a_i| \right) \left| \e\sum_{j=1}^{n} \prod_{i=1}^j |a_i^{-1}| + (w_1 -z_1) \right| .
\end{align}
The assumption is equivalent to the limit as follows 
$$ \lim_{n \rightarrow \infty} \left(\prod_{j=1}^{n} |a_j| \right)^{-\frac{1}{n}} = 1 . $$
If we put $ t_n = \prod_{j=1}^{n} |a_j^{-1}| $, then Lemma \ref{lem-key lemma1} implies that 
\begin{align} \label{eq-unboundedness of norm}
\lim_{n \rightarrow \infty}\; \cfrac{\prod_{j=1}^{n} |a_j^{-1}|}{\sum_{j=1}^{n-1} \prod_{i=1}^j |a_i^{-1}|} = 0 \quad \Longleftrightarrow \quad 
\lim_{n \rightarrow \infty} \left( \prod_{j=1}^n |a_j| \right) \sum_{j=1}^{n-1} \prod_{i=1}^j |a_i^{-1}| = \infty .
\end{align}
Observe that $ \prod_{i=1}^j |a_i^{-1}| $ is always positive and  the equation
$$ \left( \prod_{i=1}^n |a_i| \right) \sum_{j=1}^{n-1} \prod_{i=1}^j |a_i^{-1}| = \left( \prod_{i=1}^n |a_i| \right) \sum_{j=1}^{n} \prod_{i=1}^j |a_i^{-1}| -1 $$ 
is true for every $ n \in \N $. If $ \lim_{n \rightarrow \infty}\prod_{i=1}^n |a_i| < \infty $, then the lower bound of $ |w_{n+1} -z_{n+1}| $ is by \eqref{eq-inquality of norm wn-zn} as follows
$$ |w_{n+1} -z_{n+1}|
\geq \left( \prod_{i=1}^n |a_i| \right) \left| \e\sum_{j=1}^{n} \prod_{i=1}^j |a_i^{-1}| \right| - \left( \prod_{i=1}^n |a_i| \right) |w_1 -z_1|  $$
for $ n \in \N $. The first term of the right hand side is unbounded by \eqref{eq-unboundedness of norm} but the norm of the second term is bounded at any choice of $ z_1 $. In this case, the difference equation has no Hyers-Ulam stability. Observe the case that $ \lim_{n \rightarrow \infty}\prod_{i=1}^n |a_i| = \infty $. If $ \left| \e\sum_{j=1}^{n} \prod_{i=1}^j |a_i^{-1}| + (w_1 -z_1) \right| $ is unbounded or non-zero finite limit as $ n \rightarrow \infty $, then the norm $ |w_{n+1} -z_{n+1}| $ is unbounded for $ n \in \N $. The only case that we should consider is that 
$$ \lim_{n \rightarrow \infty}\prod_{i=1}^n |a_i| = \infty \ \ \text{and} \ \ \lim_{n \rightarrow \infty} \e\sum_{j=1}^{n} \prod_{i=1}^j |a_i^{-1}| + (w_1 -z_1) = 0 $$ 
for the suitable choice of $ z_1 $. The limit has the indeterminate form $ \infty \cdot 0 $. Denote $ \lim_{n \rightarrow \infty} \sum_{j=1}^{n} \prod_{i=1}^j |a_i^{-1}| $ by $  s_0 $. Choose $ z_1 = w_1 + \e s_0 $. Thus $ w_{n+1} - z_{n+1} $ is as follows
$$ w_{n+1} - z_{n+1} = \e \left( \prod_{i=1}^n a_i \right) \sum_{j=n+1}^{\infty} \prod_{i=1}^j|a_i|^{-1} . $$
It implies that 
\begin{align} \label{eq-unboundedness of wn-zn}
|w_{n+1} - z_{n+1}| &= \e \left( \prod_{i=1}^n |a_i| \right) \sum_{j=n+1}^{\infty} \prod_{i=1}^j|a_i|^{-1} .
\end{align}
Lemma \ref{lem-key lemma1 supplementary} implies that if $ \lim_{n\rightarrow \infty} t_n^{\frac{1}{n}} = 1 $, that is, $ \lim_{n\rightarrow \infty} t_n^{-\frac{1}{n}} = 1 $, then the limit $ \lim_{n\rightarrow \infty} \frac{1}{t_n}\,\sum_{j=n+1}^{\infty} t_j = \infty $ is satisfied. Put $ t_n $ to be $ \prod_{i=1}^n|a_i|^{-1} $. Then the equation \eqref{eq-unboundedness of wn-zn} is unbounded for $ n \in \N $. Hence, the absolute value of the difference, $ |w_{n+1} -z_{n+1}| $ is unbounded for $ n \in \N $ at any choice of $ z_1 $ in the complex plane. Therefore, the difference equation \eqref{eq-equation of z in sec5 thm} has no Hyers-Ulam stability. 
\end{proof}

\medskip
\noindent The limit of $ t_n^{\frac{1}{n}} $ determines Hyers-Ulam stability of the sequence generated by linear maps. The following lemma and theorem shows another result.

\medskip
\begin{lem}[Key lemma 2] \label{lem-key lemma2}
Let $ \{t_n\}_{n \in \N} $ be the sequence of positive real numbers. Then the limit $ \displaystyle \lim_{n \rightarrow \infty} t_n^{\frac{1}{n}} = 1 $ is satisfied if and only if  
$$ \lim_{n \rightarrow \infty} \frac{t_nK^n}{\sum_{j=1}^{n-1}t_jK^j} = K-1 $$
for $ K > 1 $. 
\end{lem}

\begin{proof}
The assumption is equivalent to the limit $ \lim_{n \rightarrow \infty} \frac{\log t_n}{n} = 0 $. For any given $ \delta > 0 $ we claim that there exists $ n = n(\delta) $ such that the inequalities hold as follows 
\begin{align} \label{eq-inequality of delta 2}
e^{-\delta}K \sum_{j=1}^{n-1} t_jK^j < \sum_{j=1}^{n} t_jK^j < e^{\delta}K \sum_{j=1}^{n-1} t_jK^j .
\end{align}
Suppose for contradiction that there exists $ \delta > 0 $ which satisfies that 
\begin{align} \label{eq-inequality of delta reverse 2}
 \sum_{j=1}^{n} t_jK^j \geq e^{\delta}K \sum_{j=1}^{n-1} t_jK^j \ \text{ equivalently}, \ \ t_nK^n \geq (e^{\delta}K-1)\sum_{j=1}^{n-1} t_jK^j
\end{align}
for all $ n \in \N $. We claim that 
\begin{align} \label{eq-inequality by induction 2}
t_nK^n \geq (e^{\delta}K-1)(e^{\delta}K)^{n-2}t_1K
\end{align}
for $ n = 2,3,4 \ldots $ by induction. For $ n=2 $, $ t_2K^2 \geq (e^{\delta}K-1)t_1K $ by the inequality \eqref{eq-inequality of delta reverse 2}. Assume that $ t_kK^k \geq (e^{\delta}K-1)(e^{\delta}K)^{k-2}t_1K $ for $ k = 2,3, \ldots, n $. Thus we obtain the inequality for $ t_{n+1}K^{n+1} $ as follows
\begin{align*}
t_{n+1}K^{n+1} &\geq (e^{\delta}K-1)\sum_{j=1}^{n} t_jK^j \\
&= (e^{\delta}K-1)(t_1K + t_2K^2 + \cdots + t_nK^n) \\
&\geq (e^{\delta}K-1)\big(t_1K + (e^{\delta}K-1)t_1K + (e^{\delta}K-1)(e^{\delta}K)t_1K + \\
& \qquad \cdots + (e^{\delta}K-1)(e^\delta K)^{(n-2)}t_1K \big) \\
&= (e^{\delta}K-1)\big(1 + (e^{\delta}K-1) + (e^{\delta}K-1)(e^{\delta}K) + \\
& \qquad \cdots + (e^{\delta}K-1)(e^\delta K)^{(n-2)} \big)t_1K \\
&= (e^{\delta}K-1)\big(1+ (-1+e^{\delta}K) + (-e^{\delta} K+ (e^{\delta}K)^2) + \\
& \qquad \cdots + (-(e^{\delta}K)^{n-2} + (e^{\delta}K)^{n-1}) \big)t_1 \\
&= (e^{\delta} K-1)(e^{\delta}K)^{n-1}t_1K .
\end{align*}
Thus $ t_{n+1}K^{n+1} \geq (e^{\delta}K-1)(e^{\delta}K)^{n-1}t_1K $. The claim is proved. Then the equation \eqref{eq-inequality by induction 2} implies the following inequality
\begin{align*}
\lim_{n \rightarrow \infty} \frac{\log t_n}{n} + \log K \geq \lim_{n \rightarrow \infty} \left( \frac{\log(e^{\delta K} -1)}{n} + \frac{n-2}{n} (\delta + \log K) + \frac{\log t_1K}{n} \right) .
\end{align*}
Thus the above inequality is the same as $ \log K \geq \delta + \log K $, which is the contradiction. Similarly, suppose for contradiction that 
$$ e^{-\delta}K \sum_{j=1}^{n-1} t_jK^j \geq \sum_{j=1}^{n} t_jK^j . $$
Then the contradictory inequality $ \log K \leq -\delta + \log K $ also appears. The detailed calculation is left to the readers. Then 
\eqref{eq-inequality of delta 2} implies that there exists $ n $ such that 
$$ e^{-\delta}K < \frac{\sum_{j=1}^n t_jK^j}{\sum_{j=1}^{n-1} t_jK^j} < e^{\delta}K $$
for all $ \delta > 0 $. Thus $ \delta = \delta(n) \searrow 0 $ as $ n \rightarrow \infty $. Passing the limit we obtain that $ \lim_{n \rightarrow \infty} \frac{\sum_{j=1}^n t_jK^j}{\sum_{j=1}^{n-1} t_jK^j} = K $. Then the following limit holds
$$ \lim_{n \rightarrow \infty} \frac{\sum_{j=1}^n t_jK^j}{\sum_{j=1}^{n-1} t_jK^j} - 1 
= \lim_{n \rightarrow \infty} \frac{t_nK^n}{\sum_{j=1}^{n-1}t_jK^j} = K-1 . $$
Suppose conversely that 
$$ \lim_{n \rightarrow \infty} \frac{t_nK^n}{\sum_{j=1}^{n-1}t_jK^j} = K-1 > 0 . $$
Thus for every $ n \geq N $, the following inequality holds
$$ \frac{1}{C_n}(K-1) \sum_{j=N}^{n-1}t_jK^j \leq t_nK^n \leq C_n (K-1)\sum_{j=N}^{n-1}t_jK^j $$
for some $ C_n \geq 1 $ where $ \lim_{n \rightarrow \infty} C_n = 1 $. We may assume that $ N=1 $ for simplicity. Then $ t_n $ satisfies the following inequality
\begin{align*}
t_nK^n &\leq C_n(K-1) \sum_{j=1}^{n-1}t_jK^j \\[-0.5em]
&= C_n(K-1)t_{n-1}K^{n-1} + C_n(K-1) \sum_{j=1}^{n-2}t_jK^j \\
&\leq C_n(K-1)t_{n-1}K^{n-1} + C_nC_{n-1}t_{n-1}K^{n-1} \\[1em]
&= C_n(K-1+C_{n-1})t_{n-1}K^{n-1} .
\end{align*}
Similarly, the other inequality holds for $ t_n $
$$ \frac{1}{C_n} \left( K-1+\frac{1}{C_{n-1}} \right)t_{n-1}K^{n-1} \leq t_nK^n . $$
Thus
\begin{align} \label{eq-ineq of ratio of t-n}
\frac{1}{C_nK} \left( K-1+\frac{1}{C_{n-1}} \right) \leq \frac{t_n}{t_{n-1}} \leq \frac{C_n}{K}(K-1+C_{n-1}) .
\end{align}
The equation $ \frac{t_n}{t_1} = \frac{t_2}{t_1} \frac{t_3}{t_2} \cdots \frac{t_n}{t_{n-1}} $ and the inequality \eqref{eq-ineq of ratio of t-n} implies that
$$ \left( K-1+\frac{1}{C_{j-1}} \right) \prod_{j=2}^{n} \frac{1}{C_jK} \leq \frac{t_n}{t_1} \leq (K-1+C_{j-1})\prod_{j=2}^{n} \frac{C_j}{K} . $$
%
The fact that $ \lim_{n \rightarrow \infty} C_n = 1 $, that is, $ \lim_{n \rightarrow \infty} \log C_n = 0 $ implies that $ \lim_{n \rightarrow \infty} \frac{1}{n} \sum_{j=2}^{n} \log C_j = 0 $ as Ces\`aro sum. Then $ \lim_{n \rightarrow \infty} \left[ \prod_{j=2}^{n} C_j \right]^{\frac{1}{n}} = 1 $. Using the similar calculations again $ \lim_{n \rightarrow \infty} \left[ \prod_{j=2}^{n} \left(K - 1 + C_j \right) \right]^{\frac{1}{n}} = K $. The same argument is applicable to $ \left[ \prod_{j=2}^{n}\frac{1}{C_j} \right]^{\frac{1}{n}} $ and $ \left[ \prod_{j=2}^{n} \left(K - 1 + \frac{1}{C_j}\right) \right]^{\frac{1}{n}} $. Observe that $ \lim_{n \rightarrow \infty} t_1^{\frac{1}{n}} = 1 $ for a constant $ t_1 > 0 $. Then passing the limit, we have that
$$ \lim_{n \rightarrow \infty} t_n^{\frac{1}{n}} = \frac{K}{K} =1 . $$
The proof is complete. 
\end{proof}

\medskip
\begin{lem} \label{lem-key lemma2-1}
Let $ \{t_n\}_{n \in \N} $ be the sequence of positive real numbers. Suppose that $ \displaystyle \liminf_{n \rightarrow \infty} \frac{t_n}{t_{n-1}} = 1 $ is satisfied. Then  
$$ \liminf_{n \rightarrow \infty} \frac{t_nK^n}{\sum_{j=1}^{n-1}t_jK^j} \geq K-1 $$
for $ K > 1 $. 
\end{lem}

\begin{proof}
The sequence $ \{\sum_{j=1}^{n-1}t_jK^j\}_{n \in \N} $ is increasing and its limit is $ +\infty $ as $ n \rightarrow \infty $. Thus we use Stolz-Ces\`aro thorem as follows
\begin{align} \label{eq-limit of ratio tn inf}
\liminf_{n \rightarrow \infty} \frac{\sum_{j=1}^{n}t_jK^j}{\sum_{j=1}^{n-1}t_jK^j} \geq \liminf_{n \rightarrow \infty} \frac{t_nK^n}{t_{n-1}K^{n-1}} = \liminf_{n \rightarrow \infty} \frac{t_n}{t_{n-1}} K = K .
\end{align}
Hence, the equation $ \displaystyle \frac{\sum_{j=1}^{n}t_jK^j}{\sum_{j=1}^{n-1}t_jK^j} = 1+ \frac{t_nK^n}{\sum_{j=1}^{n-1}t_jK^j} $ and \eqref{eq-limit of ratio tn inf} implies the following inequality
$$ \liminf_{n \rightarrow \infty} \frac{t_nK^n}{\sum_{j=1}^{n-1}t_jK^j} \geq K-1 . $$
\end{proof}

\begin{rk}
In Lemma \ref{lem-key lemma2-1}, the assumption $ \displaystyle \liminf_{n \rightarrow \infty} \frac{t_n}{t_{n-1}} = 1 $ cannot be relaxed to $ \displaystyle \liminf_{n \rightarrow \infty} t_n^{\frac{1}{n}} = 1 $. For instance, let $ \{t_n\}_{n \in \N} $ be the sequence defined by 
\begin{align*}
t_n = \begin{cases}
1,\quad & n= 2k-1 \\
2^k, \quad & n=2k
\end{cases}
\end{align*}
Thus $ \displaystyle \liminf_{n \rightarrow \infty} t_n^{\frac{1}{n}} = 1 $. However, 
\begin{align*}
\liminf_{n \rightarrow \infty} \frac{t_nK^n}{\sum\limits_{j=1}^{n-1}t_jK^j} &= \lim_{k \rightarrow \infty} \min \left\{ \frac{2^kK^{2k}}{\sum\limits_{j=1}^k K^{2j-1} + \sum\limits_{j=1}^k 2^jK^{2j}} ,\ \ \frac{K^{2k+1}}{\sum\limits_{j=1}^{k+1} K^{2j-1} + \sum\limits_{j=1}^k 2^jK^{2j}}  \right\} \\[0.8em]
&= \lim_{k \rightarrow \infty} \frac{K^{2k+1}}{\dfrac{K(K^{2k}-1)}{K^2-1} + \dfrac{2K^2(2^kK^{2k}-1)}{2K^2-1}} \\[0.3em]
&= 0 \neq K-1 .
\end{align*}
\end{rk}

\section{Exponential growth rate and HU stability}
Exponential growth rate of linear recurrence is expressed in terms of the limit related to $ a_n $ of the map $ g_n(z) = a_nz+ b_n $. In particular, If the limit of $ \left(\prod_{j=1}^{n-1} |a_j| \right)^{\frac{1}{n}} $ is greater or less than one, then we show that the difference equation $ z_{n+1} = g_n(z_n) $ has Hyers-Ulam stability. 
\medskip
\begin{thm} \label{thm-stability of exponential shrinking}
Let $ \{z_n\}_{n\in \N} $ be the sequence defined as follows
\begin{align} \label{eq-equation of z in sec6}
z_{n+1} = a_nz_n + b_n
\end{align}
for $ n \in \N $. Suppose that 
$$ \lim_{n \rightarrow \infty} \left(\prod_{j=1}^{n} |a_j| \right) ^{\frac{1}{n}} = \frac{1}{K} < 1 . $$
Then the difference equation \eqref{eq-equation of z in sec6} has Hyers-Ulam stability. 
\end{thm}

\begin{proof}
The assumption is equivalent to the following limit 
$$ \lim_{n \rightarrow \infty} \left(\prod_{j=1}^{n} \left| a_j^{-1} \right| \right)^{\frac{1}{n}} = {K} > 1 . $$ Thus we may assume that 
$$ \prod_{j=1}^{n} |a_j^{-1}| = t_nK^n \ \ \text{and} \ \ \lim_{n \rightarrow \infty} t_n^{\frac{1}{n}} = 1 $$
where $ \{ t_n \}_{n \in \N} $ is the sequence of positive real numbers. Lemma \ref{lem-key lemma2} implies that 
\begin{align} \label{eq-equiv limit 1}
& \lim_{n \rightarrow \infty}\; \cfrac{\prod_{j=1}^{n} |a_j^{-1}|}{\sum\limits_{j=1}^{n-1} \prod_{i=1}^{j} |a_i^{-1}|} = K - 1 > 0  
\qquad \text{if and only if} \nonumber \\[0.3em]
& \lim_{n \rightarrow \infty} \left( \prod_{j=1}^{n} |a_j| \right) \sum_{j=1}^{n-1} \prod_{i=1}^{j} |a_i^{-1}| = \lim_{n \rightarrow \infty} \sum_{j=1}^{n-1} \prod_{i=j+1}^{n} |a_i| = 
\frac{1}{K-1} > 0 .
\end{align}
Define the map $ g_n $ as $ g_n(z) = a_nz + b_n $ for $ n \in \N $. Let $ \{ w_n \}_{n \in \N} $ be the sequence which satisfies the equations $ w_{n+1} = g_n(w_n) + r_n $ for $ n \in \N $ where $ | r_n | \leq \e $ for all $ n \in \N $. Proposition \ref{prop-stability of linear maps} and the equation \eqref{eq-difference between w and z again} implies that 
\begin{align*} 
w_{n+1} -z_{n+1} &= \sum_{j=0}^{n-1} r_{n-j} \prod_{i=1}^j a_{n-j+i} + (w_1 -z_1) \prod_{j=1}^n a_j \\
&= \sum_{j=1}^{n} r_{j} \prod_{i=j+1}^n a_{i} + (w_1 -z_1) \prod_{j=1}^n a_j . 
\end{align*}
Choose $ z_1 = w_1 $. Then the equation \eqref{eq-equiv limit 1} yields that
\begin{align*}
|w_{n+1} -z_{n+1}| &\leq \left| \sum_{j=1}^n r_j \prod_{i=j+1}^na_i \right| \\
&\leq \e \sum_{j=1}^n  \prod_{i=j+1}^n | a_i | \\
&\leq \e \left( \sum_{j=1}^{n-1} \prod_{i=j+1}^n | a_i | +1 \right) \\
&\leq \e \left( \frac{1}{K-1} + 1\right) \\
&\leq \frac{K\e}{K-1}
\end{align*}
for all $ n \in \N_0 $. Hence, the difference equation \eqref{eq-equation of z in sec6} has Hyers-Ulam stability. 
\end{proof}

\medskip
\begin{thm} \label{thm-stability of exponential expanding}
Let $ \{z_n\}_{n\in \N} $ be the sequence defined as follows
\begin{align} \label{eq-equation of z2 in sec6}
z_{n+1} = a_nz_n + b_n
\end{align}
for $ n \in \N $. Suppose that 
$$ \lim_{n \rightarrow \infty} \left(\prod_{j=1}^{n} |a_j| \right)^{\frac{1}{n}} = {K} > 1 . $$
Then the difference equation \eqref{eq-equation of z2 in sec6} has Hyers-Ulam stability. 
\end{thm}

\begin{proof}
we may assume that 
$$ \prod_{j=1}^{n} |a_j| = t_nK^n \ \ \text{and} \ \ \lim_{n \rightarrow \infty} t_n^{\frac{1}{n}} = 1 $$
where $ \{ t_n \}_{n \in \N} $ is the sequence of positive real numbers. 
%
Proposition \ref{prop-stability of linear maps} and the equation \eqref{eq-difference between w and z again} implies that 
\begin{align} \label{eq-difference between w and z sec6}
w_{n+1} -z_{n+1} &= \sum_{j=0}^{n-1} r_{n-j} \prod_{i=1}^j a_{n-j+i} + (w_1 -z_1) \prod_{j=1}^n a_j \nonumber \\
&= \sum_{j=1}^{n} r_{j} \prod_{i=j+1}^n a_{i} + (w_1 -z_1) \prod_{j=1}^n a_j \nonumber \\
&= \left(\prod_{j=1}^n a_j \right) \left( \sum_{j=1}^{n} r_{j} \prod_{i=1}^j a_{i}^{-1} + w_1 -z_1 \right). 
\end{align}
Let $ \{w_n\}_{n\in \N} $ be the sequence defined by the recursion formula 
$$ w_{n+1} = a_nw_n + b_n + r_n $$
where $ |r_n| \leq \e $ for all $ n \in \N $. For a given sequence $ \{w_n\}_{n\in \N} $, define $ z_1 $ as follows 
\begin{align} \label{eq-choice of the initial term K big}
z_1 = w_1 + \sum_{j=1}^{\infty} r_{j} \prod_{i=1}^j a_{i}^{-1} .
\end{align}
%
Since the sequence $ \sum_{j=1}^{\infty} \prod_{i=1}^j |a_{i}^{-1}| $ is convergent by the root test, the number $ z_1 $ is well defined. Apply the defined $ z_1 $ to the equation \eqref{eq-difference between w and z sec6} as follows
\begin{align} 
|w_{n+1} -z_{n+1}| &= \left(\prod_{j=1}^n |a_j| \right) \left| \sum_{j=1}^{n} r_{j} \prod_{i=1}^j a_{i}^{-1} - \sum_{j=1}^{\infty} r_{j} \prod_{i=1}^j a_{i}^{-1} \right| \nonumber \\
&= \left(\prod_{j=1}^n |a_j| \right) \left| - \sum_{j=n+1}^{\infty} r_{j} \prod_{i=1}^j a_{i}^{-1} \right| \nonumber \\
&\leq \e \left(\prod_{j=1}^n |a_j| \right) \sum_{j=n+1}^{\infty} \prod_{i=1}^{j} |a_{i}^{-1}| = \e \sum_{j=n+1}^{\infty} \prod_{i=n+1}^{j} |a_{i}^{-1}| \label{eq-bounds of diff of sequence2} \\
& = \e \left(t_nK^n \right) \sum_{j=n+1}^{\infty} \frac{1}{t_j K^j} .\label{eq-bounds of diff of sequences}
\end{align}
Let $ L $ is a positive number between $ 1 $ and $ K $, that is, $ K > L > 1 $. Denote $ t_n\left( \frac{L}{K} \right)^{n} $ by $ s_n $ for $ n \in \N $. Since the inequality $ \dfrac{1}{K^{j-n}} <  \dfrac{1}{L^{j-n}} $ is satisfied for all $ j > n $, the following inequality holds
\begin{align*}
\sum_{j=n+1}^{\infty}\frac{t_{n}K^{n}}{t_j K^j} &< \sum_{j=n+1}^{\infty}\frac{t_{n}L^{n}}{t_j L^j} \\
&= t_{n}\left( \frac{L}{K} \right)^{n}K^{n} \sum_{j=n+1}^{\infty} \frac{1}{t_j \left( \frac{L}{K} \right)^{j}K^j} \\
&= s_{n}K^{n} \sum_{j=n+1}^{\infty} \frac{1}{s_j K^j} .
\end{align*}
Observe that $ \lim_{n \rightarrow \infty} s_n  = 0 $ and $ \lim_{n \rightarrow \infty} (s_n)^{\frac{1}{n}}  = \frac{L}{K} < 1 $. 
In order to estimate $ \e \left(t_{n}K^{n} \right) \sum_{j=n+1}^{\infty} \frac{1}{t_j K^j} $, we define $ C^1 $ function on the set of positive real numbers, say $ g $, which satisfies the following properties.
\begin{enumerate}
\item $ g(x) > 0 $ for all $ x > 0 $,
\item $ g(j) = \frac{1}{s_j} $ for every $ j \in \N $,
\item $ \lim_{n \rightarrow \infty} g(x) = \infty $, and 
\item $ \lim_{x \rightarrow \infty} \{g(x)\}^{\frac{1}{x}} = \frac{K}{L} > 1 $.
\end{enumerate}
The forth property is equivalent to the limit $ \lim_{x \rightarrow \infty} \frac{\ln(g(x))}{x} = \ln \left(\frac{K}{L} \right) $ or 
\begin{align} \label{eq-limit of g}
\lim_{x \rightarrow \infty} \frac{g'(x)}{g(x)} = \ln {K} - \ln {L} > 0
\end{align}
by L'Hospital's rule. Consider the function $ g(x)K^{-x} $ and calculate its derivative 
$$ \left( g(x)K^{-x}\right)' = \left[ \frac{g'(x)}{g(x)} - \ln K \right] g(x)K^{-x} . $$
The equation \eqref{eq-limit of g} implies that there exists big enough $ x_0 $ such that the inequality $ \frac{|g'(x)|}{g(x)} - \ln K < 0 $ for all $ x > x_0 $. Thus $ g $ is decreasing on the interval $ \{ x\ | \ x > x_0 \} $. Moreover, there exists a positive integer $ n $ as follows
\begin{align*}
\frac{|g'(x)|}{g(x)} \leq \frac{1}{2}\ln K 
\end{align*}
for all $ x \geq n \geq x_0 $ with a suitable choice of $ L > \sqrt{K} $. 
Then we may use the comparison between infinite series and definite integral as follows 
\begin{align} \label{eq-comparison between series and integral}
\sum_{j=n+1}^{\infty} \frac{1}{s_j K^j} &\leq \int_{n}^{\infty}g(x)K^{-x}dx .
\end{align}
The definite integral on $ [n-1,\infty ) $ is estimated as follows 
\begin{align} \label{eq-bounds of integral}
\int_{n}^{\infty}g(x)K^{-x}dx 
&= \int_{n}^{\infty}g(x)e^{-x\ln K}dx \nonumber \\[0.5em]
&= \left[ -\frac{1}{\ln K}\, g(x)e^{-x\ln K} \right]_{n}^{\infty} + \frac{1}{\ln K} \int_{n}^{\infty}g'(x)e^{-x\ln K}dx \nonumber \\[0.5em]
&\leq \frac{1}{\ln K}\, g(n)e^{-n\ln K} + \frac{1}{2}\int_{n}^{\infty}g(x)K^{-x}dx .
\end{align}
Thus an upper bound of the integral \eqref{eq-bounds of integral} is as follows
\begin{align*}
\int_{n}^{\infty}g(x)K^{-x}dx \leq \frac{2}{\ln K}\,\frac{1}{s_{n}}K^{-n} .
\end{align*}
The estimations \eqref{eq-comparison between series and integral}  are applied to the difference $ |w_{n+1} -z_{n+1}| $ in the equation \eqref{eq-bounds of diff of sequences}. Then we have that
\begin{align*}
|w_{n+1} -z_{n+1}| &\leq \e \left(t_{n}K^{n} \right) \sum_{j=n+1}^{\infty} \frac{1}{t_j K^j} \\[0.2em]
&< \e \left(s_{n}K^{n} \right) \sum_{j=n+1}^{\infty} \frac{1}{s_j K^j} \\[0.2em]
&\leq \e \left(s_{n}K^{n} \right) \int_{n}^{\infty}g(x)K^{-x}dx \\[0.5em]
&= \e \left(s_{n}K^{n} \right) \frac{2}{\ln K}\,\frac{1}{s_{n}}K^{-n} \\
&= \frac{2\e}{\ln K} ,
\end{align*}
that is, we obtain an upper bound of difference between $ w_{n+1} $ and $ z_{n+1} $ as follows 
\begin{align} \label{eq-upper bound for HU stability}
|w_{n+1} -z_{n+1}| \leq \frac{2\e}{\ln K}
\end{align}
for all sufficiently large $ n \in \N $. An upper bound of $ |w_N -z_N | $ is estimated for $ N=1,2,\ldots,n $ by the equation \eqref{eq-bounds of diff of sequence2} and the estimation \eqref{eq-upper bound for HU stability} as follows
\begin{align*}
|w_N-z_N| &\leq \e \sum_{j=N}^{\infty} \prod_{i=N}^{j} |a_{i}^{-1}| \\
&= \e \left(\;\sum_{j=N}^{n} \prod_{i=N}^{j} |a_{i}^{-1}| + \sum_{j=n+1}^{\infty} \prod_{i=N}^{j} |a_{i}^{-1}| \right) \\[0.3em]
&= \e \left[\; \sum_{j=N}^{n} \prod_{i=N}^{j} |a_{i}^{-1}| + \left(\prod_{i=N}^{n} |a_{i}^{-1}| \right)\sum_{j=n+1}^{\infty} \prod_{i=n+1}^{j} |a_{i}^{-1}| \; \right] \\[0.3em]
&\leq \e \left[\; \sum_{j=N}^{n} \prod_{i=N}^{j} |a_{i}^{-1}| + \left(\prod_{i=N}^{n} |a_{i}^{-1}| \right) \frac{2}{\ln K}\right]
\end{align*} 
for $ N =1,2,\ldots,n $. Hence, the difference equation \eqref{eq-equation of z2 in sec6} has Hyers-Ulam stability. 
\end{proof}

\medskip
\noindent The following theorem is a counterpart of Theorem \ref{thm-HU stability with periodic coefficients 1}, which is Hyers-Ulam stability of the sequence of generated by linear maps with periodic coefficients.

\medskip

\begin{thm} \label{thm-HU stability with periodic coefficients 2}
Let $ \{z_n\}_{n\in \N} $ be the sequence defined as follows
\begin{align} \label{eq-equation of z3 in sec6}
z_{n+1} = a_nz_n + b_n
\end{align}
for $ n \in \N $. Let $ a_n $ be the periodic coefficient, that is, $ a_{n+p} = a_n $ for some $ p \in \N $ and for every $ n \in \N $. Assume that $ \prod_{k=1}^p \left| a_k \right| > 1 $. Then the difference equation \eqref{eq-equation of z3 in sec6} has Hyers-Ulam stability.
\end{thm}

\begin{proof}
It suffice to show that 
$$ \lim_{n \rightarrow \infty} \left(\prod_{j=1}^{n} |a_j| \right)^{\frac{1}{n}} = {K} $$
for some $ K > 1 $ by Theorem \ref{thm-stability of exponential expanding}. Denote $ \prod_{k=1}^p |a_k| $ by $ K^p $ for some $ K > 1 $. 
Thus
\begin{align} \label{eq-periodic coefficient expanding}
\left(\prod_{j=1}^{n} |a_j| \right)^{\frac{1}{n}} &= \left(\prod_{j=1}^{mp+l}|a_j| \right)^{\frac{1}{mp+l}} \nonumber \\
&=  \left[ K^{mp}|a_1a_2\cdots a_l| \right]^{\frac{1}{mp+l}} \nonumber \\
&= K \left[ K^{-l}|a_1a_2\cdots a_l| \right]^{\frac{1}{mp+l}}
\end{align}
where $ n = mp + l $ for $ l =0,1,2,\ldots,p-1 $. Then the limit of the equation \eqref{eq-periodic coefficient expanding} is
$$ \lim_{m \rightarrow \infty} K \left[ K^{-l}|a_1a_2\cdots a_l| \right]^{\frac{1}{mp+l}} = K $$
for each $ l =0,1,2,\ldots,p-1 $. Hence, the difference equation \eqref{eq-equation of z3 in sec6} has Hyers-Ulam stability by Theorem \ref{thm-stability of exponential expanding}.
\end{proof}

\medskip
\begin{exa}
Let $ \{u_n\}_{n \in \N} $ be the sequence which satisfies the equation
\begin{align} \label{eq-example1 in sec6}
u_{n+1} &= a_nu_n + 5 \quad \text{where} \quad a_n =
\begin{cases}
3, \quad n = pk  \\
1, \quad \text{otherwise}
\end{cases}
\end{align}
for some $ p \in \N $ and for every $ k \in \N $. Thus $ \prod_{j=1}^p |a_j| = 3 > 1 $. Theorem \ref{thm-HU stability with periodic coefficients 2} implies that the difference equation \eqref{eq-example1 in sec6} has Hyers-Ulam stability. \\
Let $ \{v_n\}_{n \in \N} $ be the sequence which satisfies the equation
\begin{align} \label{eq-example2 in sec6}
v_{n+1} &= a_nv_n + 5 \quad \text{where} \quad a_n =
\begin{cases}
3, \quad n = k^2  \\
1, \quad \text{otherwise}
\end{cases}
\end{align}
for every $ k \in \N $. Thus $ \prod_{j=1}^{n^2} |a_j| = 3^n $. Then 
$$ \lim_{n \rightarrow \infty} \left(\prod_{j=1}^{n^2} |a_j|\right)^{\frac{1}{n^2}} = \lim_{n \rightarrow \infty} 3^{\frac{1}{n}} = 1 . $$ 
Theorem \ref{thm-no HU stability subexp} implies that the difference equation \eqref{eq-example2 in sec6} has no Hyers-Ulam stability. 
\end{exa}

\medskip

\begin{exa}
Let $ \{z_n\}_{n \in \N} $ be the sequence which satisfies the equation
\begin{align} \label{eq-example2 in sec6}
z_{n+1} &= a_nz_n + 3 \quad \text{where} \quad a_n = K \left(\frac{n}{n-1}\right)^{\alpha}
\end{align}
for $ n =2,3,4,\ldots $ and $ a_1 = 3 $ where $ \alpha > 0 $ and $ K>1 $. Thus we obtain $ \prod_{j=1}^{n} |a_j| = n^{\alpha}K^n $. Observe that $ \lim_{n \rightarrow \infty}\left(\prod_{j=1}^{n} |a_j| \right)^{\frac{1}{n}} = K > 1 $. We may choose $ t_n = n^{\alpha} $ for every $ n \in \N $. Then we follow the proof of Theorem \ref{thm-stability of exponential expanding} to show Hyers-Ulam stability of the difference equation \eqref{eq-example2 in sec6}. Let $ g(x) = \frac{1}{x^{\alpha}} $. Thus $ g(j) = \frac{1}{t_j} $ for every $ j \in \N $. Moreover, $ g(x)K^{-x} $ is a decreasing function on $ x > 0 $ because
$$ \left(\frac{K^{-x}}{x^{\alpha}}\right)' = -\frac{K^{-x}}{x^{\alpha}}\left( \frac{\alpha}{x}+\ln K \right) < 0 $$ 
for all $ x > 0 $. Let us calculate an upper bound of the improper integral of the function as follows
\begin{align*}
\int_n^{\infty} \frac{K^{-x}}{x^{\alpha}}\,dx &= \int_n^{\infty} \frac{1}{x^{\alpha}}\,e^{-x\ln K}dx \\[0.3em]
&= \left[\;-\frac{1}{\ln K}\cdot\frac{e^{-x\ln K}}{x^{\alpha}}\; \right] - \int_n^{\infty} \frac{\alpha}{\ln K}\cdot\frac{e^{-x\ln K}}{x^{{\alpha}+1}}\,dx \\[0.3em]
&< \frac{1}{\ln K}\cdot\frac{e^{-n\ln K}}{n^{\alpha}} \\[0.3em]
&= \frac{1}{\ln K}\cdot\frac{K^{-n}}{n^{\alpha}}
\end{align*}
for every $ n \in \N $. For a given $ \e > 0 $. if the given sequence $ \{w_n\}_{n \in \N} $ satisfies that
$$ |w_{n+1} - a_nw_n - 3| \leq \e
$$
for all $ n \in \N $, then the suitable choice of $ z_1 $ implies the following bounds by the equation \eqref{eq-bounds of diff of sequences} and \eqref{eq-comparison between series and integral} 
\begin{align*}
|w_{n+1}-z_{n+1}| \leq \e \left(t_nK^n \right) \sum_{j=n+1}^{\infty} \frac{1}{t_j K^j} \leq \e n^{\alpha}K^n \cdot \frac{1}{\ln K}\cdot\frac{K^{-n}}{n^{\alpha}} = \frac{\e}{\ln K}
\end{align*}
for all $ n \in \N $. The equation \eqref{eq-choice of the initial term K big} implies that
$$ |w_1 - z_1| \leq \e\sum_{j=1}^{\infty}\prod_{i=1}^j |a_i^{-1}| = \e\frac{1}{|a_1|}\sum_{j=2}^{\infty}\prod_{i=1}^j |a_i^{-1}|= \frac{1}{3}\cdot\frac{\e}{\ln K}
$$
Hence, we conclude that 
$$ |w_{n}-z_{n}| \leq \frac{\e}{\ln K} $$
for all $ n \in \N $. The difference equation \eqref{eq-example2 in sec6} has Hyers-Ulam stability. 
\end{exa}

\section{Application to the higher order equation}
In this section we suggest the sufficient condition for Hyers-Ulam stability of the second order difference equation of a particular type with periodic coefficients. 
\begin{exa}
Let $ \{z_n\}_{n \in \N} $ be the sequence which satisfies the equation
\begin{align}
z_{n+1} &= a_nz_n + b_n \label{eq-first order difference eq} \\ 
&\text{where} \quad a_n =
\begin{cases}
s, \quad n = 2k + 1 \\
q, \quad n = 2k 
\end{cases}
\ \text{and} \quad b_n =
\begin{cases}
u, \quad n = 2k + 1 \\
v, \quad n = 2k 
\end{cases} \nonumber
\end{align}
where $ s,q,u $ and $ v $ are complex numbers for $ n \in \N\cup\{0\} $. Suppose that $ |sq| \neq 1 $, that is, $ sq $ is not contained in the unit circle in complex plane. Then Theorem \ref{thm-HU stability with periodic coefficients 1} or Theorem \ref{thm-HU stability with periodic coefficients 2} implies that $ \{z_n\}_{n \in \N} $ has Hyers-Ulam stability. Thus for a given $ \e > 0 $, the sequence $ \{w_n\}_{n \in \N} $ satisfies that 
\begin{align} \label{eq-first order difference ineq}
w_{n+1} = a_nw_n + b_n + r_n
\end{align}
where $ |r_n| \leq \e $ for all $ n \in \N\cup\{0,-1\} $. Then $ |w_n -z_n| = G(\e) $ is bounded and $ G(\e) \rightarrow 0 $ as $ \e \rightarrow 0 $. However, the equation \eqref{eq-first order difference eq} implies the second order difference equation $ z_{n+1} = (a_n-1)z_n + a_{n-1}z_{n-1} + u+v $, that is, 
\begin{align} \label{eq-second order difference particular type}
\begin{cases}
z_{2k+1} = (q-1)z_{2k} + sz_{2k-1} + v+u \\
z_{2k+2} = (s-1)z_{2k+1} + qz_{2k} + u+v 
\end{cases}
\end{align}
for $ k \in \N\cup\{0\} $ with the initial values $ z_{-1} $ and $ z_0 $. Define $ a_{-1} = s $ and $ b_{-1} $ is the numbers which satisfy the equation $ z_0 = a_{-1}z_{-1}+b_{-1} $. Then the equation \eqref{eq-first order difference ineq} satisfies the following second order difference equation
$$ w_{n+1} = (a_n-1)w_n + a_{n-1}w_{n-1} + u+v + r_n +r_{n-1} $$
for $ n \in \N\cup\{0\} $. Since $ |r_n +r_{n-1}| \leq 2\e $, the equation \eqref{eq-second order difference particular type} has Hyers-Ulam stability also. 
\end{exa}

\medskip

\section{Conclusion and final note}
{\em Subexponential growth rate} of the sequence of positive real numbers is suggested. For the linear difference equation $ z_{n+1} = a_nz_n + b_n $, we use a sufficient condition $ \lim_{n\rightarrow \infty} \prod_{j=1}^n |a_j|^{\frac{1}{n}} =1 $ to show that the linear recurrence does not have Hyers-Ulam stability. Moreover, if the quantity, $ \lim_{n\rightarrow \infty} \prod_{j=1}^n |a_j|^{\frac{1}{n}} $ is greater than one or less than one, then $ |z_n| $ eventually expands or shrinks at least exponentially fast for $ n \in \N $. Thus Lemma \ref{lem-key lemma2} is essential to prove that linear recurrence is Hyers-Ulam stable. The existence of the limit, $ \lim_{n\rightarrow \infty} \prod_{j=1}^n |a_j|^{\frac{1}{n}} $ provide the unified approach for us to determine Hyers-Ulam stability of one dimensional linear recurrence. \\
\noindent There are something to discuss. 
We do not know yet whether the subexponential growth of linear recurrence requires the convergence of $ \prod_{j=1}^n |a_j|^{\frac{1}{n}} $ as $ n \rightarrow \infty $ or not. Thus the converse of Lemma \ref{lem-key lemma1} is neither proved nor disproved. 
In higher dimension, the subexponential growth of the sequence $ \{ \| t_n\| \}_{n\in \N} $ can be defined as follows
$$ \lim_{n\rightarrow \infty} \frac{\| t_n \|}{\left\| \sum_{j=1}^{n-1}t_j \right\|} =0 \quad \text{or} \ \ \lim_{n\rightarrow \infty} \frac{\| t_n \|}{ \sum_{j=1}^{n-1} \|t_j\| } =0 $$
where $ t_n $ is the linear operator from the higher dimensional space to itself for $ n \in \N $. 



\end{document}